\input ./style/arxiv-general.cfg
\documentclass[aos,MSNbibl,seceqn,dvips]{arximspdf}
\makeatletter
   \@ifpackageloaded{graphicx}{}{\usepackage{graphicx}}
\makeatother

\doi{10.1214/15-AOS1343}
\volume{44}
\issue{1}
\pubyear{2016}
\firstpage{31}
\lastpage{57}
\docsubty{FLA}

\makeatletter
\newcommand{\rrvert}{\vert}
\newcommand{\rrVert}{\Vert}
\newcommand{\llvert}{\vert}
\newcommand{\llVert}{\Vert}
\newcommand{\eqref}[1]{(\ref{#1})}
\newtheorem{theorem}{Theorem} 
\newtheorem{lemma}{Lemma} 
\newtheorem{proposition}{Proposition} 
\newproclaim{remark}{Remark} 
\newcommand{\E}{\mathbb{E}}
\newcommand{\R}{\mathbb{R}}
\newcommand{\N}{\mathbb{N}}
\renewcommand{\P}{\mathbb{P}}
\newcommand{\M}{\mathcal{M}}
\makeatother

\begin{document}
\begin{frontmatter}

\title{Asymptotics in directed exponential random graph models with an
increasing bi-degree sequence}
\runtitle{Directed exponential random graph models}

\begin{aug}
\author[A]{\fnms{Ting}~\snm{Yan}\thanksref{T1}\ead[label=e1]{tingyanty@mail.ccnu.edu.cn}},
\author[B]{\fnms{Chenlei}~\snm{Leng}\ead[label=e2]{C.leng@warwick.ac.uk}}
\and
\author[C]{\fnms{Ji}~\snm{Zhu}\corref{}\thanksref{T3}\ead[label=e3]{jizhu@umich.edu}\ead[label=u1,url]{http://www.foo.com}}
\runauthor{T. Yan, C. Leng and J. Zhu}
\affiliation{Central China Normal University,
University of Warwick and \\
University of Michigan}
\address[A]{T. Yan\\
Department of Statistics\\
Central China Normal University \\
Wuhan, 430079\\
China\\
\printead{e1}}
\address[B]{C. Leng\\
Department of Statistics\\
University of Warwick\\
Coventry, CV4 7AL\\
United Kingdom\\
\printead{e2}}
\address[C]{J. Zhu\\
Department of Statistics\\
University of Michigan\\
439 West Hall\\
1085 S. University Ave\\
Ann Arbor, Michigan 48109-1107\\
USA\\
\printead{e3}}
\end{aug}
\thankstext{T1}{Supported in part by the National Science Foundation
of China (No. 11401239).}
\thankstext{T3}{Supported in part by NSF Grant DMS-14-07698.}

%
\received{\smonth{12} \syear{2014}}
%
\revised{\smonth{5} \syear{2015}}

%
\begin{abstract}
Although asymptotic analyses of undirected network models based on
degree sequences have started to
appear in recent literature, it remains an open problem to study
statistical properties of
directed network models. In this paper, we provide for the first time a
rigorous analysis
of directed exponential random graph models using the in-degrees and
out-degrees as
sufficient statistics with binary as well as continuous weighted edges.
We establish
the uniform consistency and the asymptotic normality for the maximum likelihood
estimate, when the number of parameters grows and only one realized observation
of the graph is available. One key technique in the proofs is to approximate
the inverse of the Fisher information matrix using a simple matrix with
high accuracy.
Numerical studies confirm our theoretical findings.
\end{abstract}

%
\end{frontmatter}

\section{Introduction}
Recent advances in computing and measurement technologies have led to
an explosion in the amount of data with network structures in a variety
of fields including social networks
\cite{GirvanNewman2002,KossinetsWatts2006}, communication networks
\cite{AdamicGlance2005,DiesnerCarley2005,Akoglu.et.al.2012},
biological networks \cite
{vonMering.et.al.2002,BaderHogue2003,NepuszYuPaccanaro2012},
disease transmission networks \cite{Newman2002,Salathea2010} and so
on. This creates an urgent need to understand the generative mechanism
of these networks and to explore various characteristics of the network
structures in a principled way. Statistical models are useful tools to
this end, since they can capture the regularities of network processes
and variability of network configurations of interests,
and help to understand the uncertainty associated with observed
outcomes \cite{Robins.et.al.2007a,Robins.et.al.2007b}.
At the same time, data with network structures pose new challenges for
statistical inference, in particular asymptotic analysis when
only one realized network is observed and one is often interested in
the asymptotic phenomena with the growing size of the network \cite
{Fienberg2012}.

The in- and out-degrees of vertices (or degrees for undirected
networks) preliminarily summarize the information contained in a
network, and their distributions provide important insights for
understanding the generative mechanism of networks. In the undirected
case, the degree sequence has been extensively studied \cite
{ChatterjeeDiaconisSly2011,BickelChenLevina2012,OlhedeWolfe2012,ZhaoLevinaZhu2012,Rinaldo2013,HillarWibisono2013}.
In particular, its distributions have been explored under the framework
of the exponential family parameterized by the so-called ``potentials''
of vertices recently, for example, the ``$\beta$-model'' by \cite
{ChatterjeeDiaconisSly2011} for binary edges or ``maximum entropy
models'' by \cite{HillarWibisono2013} for weighted edges in which
the degree sequence is the exclusively sufficient statistic.
It is also worth to note that the asymptotic theory of the maximum
likelihood estimates (MLEs) for these models have not been derived
until very recently \cite
{ChatterjeeDiaconisSly2011,HillarWibisono2013,YanXu2013,YanZhaoQin2013}.
In the directed case,
how to construct and sample directed graphs with given in- and
out-degree (sometimes referred as ``bi-degree'') sequences have been
studied \cite
{ErdosPeterMiklosToroczkai2010,ChenOlvera-Cravioto2013,KimGenio2012}.
However, statistical inference is not available, especially for
asymptotic analysis. 
The distributions of the bi-degrees were studied in \cite
{Robins.et.al.2009} through empirical examples for social networks, but
the work lacked theoretical analysis. 

In this paper, we study the distribution of the bi-degree sequence when
it is the sufficient statistic in a directed graph.
Recall the Koopman--Pitman--Darmois theorem or the principle of maximum
entropy \cite{WainwrightJordan2008,Wu1997}, which states that the
probability mass function of the bi-degree sequence must admit the form
of the exponential family. We will characterize the exponential family
distributions for the bi-degree sequence with three types of weighted
edges (binary, discrete and continuous) and conduct the maximum
likelihood inference.

In the model we study, one out-degree parameter and one in-degree
parameter are needed for each vertex. As a result, the total number of
parameters is twice of the number of the vertices. As the size of the
network increases, the number of parameters goes to infinity.
This makes asymptotic inference very challenging.
Establishing the uniform consistency and asymptotic normality of the
MLE are the aims of this paper. To the best of our knowledge, it is the
first time that such results are derived in directed exponential random
graph models with weighted edges. We remark further that our proofs are
highly nontrivial. One key feature of our proofs lies in approximating
the inverse of the Fisher information matrix by a simple matrix with
small approximation errors. This approximation is utilized to derive a
Newton iterative algorithm with geometrically fast rate of convergence,
which leads to the proof of uniform consistency, and it is also
utilized to derive approximately explicit expressions of the
estimators, which leads to the proof of asymptotic normality.
Furthermore, the approximate inverse makes the asymptotic variances of
estimators explicit and concise.
We note that \cite{Haberman1977,Haberman1981} have studied problems
related to the present paper but the methods therein cannot be applied
to the model we study. This is explained in detail at the end of the
next section after we state the main theorems.

Next, we formally describe the models considered in this paper. 
Consider a directed graph $\mathcal{G}$ on $n \geq2$ vertices labeled
by $1, \ldots, n$. Let $a_{i,j}\in\Omega$ be the weight of the
directed edge from $i$ to $j$, where $\Omega\subseteq\R$ is the set
of all possible weight values, and $A=(a_{i,j})$ be the adjacency
matrix of $\mathcal{G}$. We consider three cases: $\Omega=\{0,1\}$,
$\Omega=[0, \infty)$ and $\Omega=\{0, 1, 2, \ldots\}$, where the
first case is the usual binary edge.
We assume that there are no self-loops, that is, $a_{i,i}=0$. Let
$d_{i}= \sum_{j \neq i} a_{i,j}$ be the out-degree of vertex $i$ and
$\mathbf{d}=(d_1, \ldots, d_n)^\top$ be the out-degree sequence of
the graph $\mathcal{G}$. Similarly, define $b_j = \sum_{i \neq j}
a_{i,j}$ as the in-degree of vertex $j$ and $\mathbf{b}=(b_1, \ldots,
b_n)^\top$ as the in-degree sequence. The pair $\{\mathbf{b}, \mathbf
{d}\}$ or $\{(b_1, d_1), \ldots, (b_n, d_n)\}$ are the bi-degree sequence.
Then the density or probability mass function on $\mathcal{G}$
parameterized by exponential family distributions with respect to some
canonical measure $\nu$ is
%
\begin{equation}
\label{Eq:density:whole} p(\mathcal{G}) = \exp\bigl( \bolds{\alpha
}^\top
\mathbf{d} + \bolds{\beta}^\top\mathbf{b} - Z(\bolds{\alpha}, \bolds{
\beta}) \bigr),
\end{equation}
where $Z(\bolds{\alpha}, \bolds{\beta})$ is the log-partition function,
$\bolds{\alpha}=(\alpha_1, \ldots, \alpha_n)^\top$ is a parameter
vector tied to the out-degree sequence, and $\bolds{\beta}=(\beta_1,
\ldots, \beta_n)^\top$ is a parameter vector tied to the in-degree sequence.
This model can be viewed 
as a directed version of the $\beta$-model \cite
{ChatterjeeDiaconisSly2011}.
It can also be represented as the log-linear model \cite
{FienbergWasserman1981a,FienbergWasserman1981b,FienbergRinaldo2012}
and the algorithm developed for the log-linear model can be used to
compute the MLE.
As explained by \cite{HollandLeinhardt1981}, $\alpha_i$ quantifies
the effect of an outgoing edge from vertex $i$ and $\beta_j$
quantifies the effect of an incoming edge connecting to vertex $j$. If
$\alpha_i$ is large and positive, vertex $i$ will tend to have a
relatively large out-degree. Similarly, if $\beta_j$ is large and
positive, vertex $j$ tends to have a relatively large in-degree. Note that
%
\begin{eqnarray}\label{Eq:tran}
\exp\bigl( \bolds{\alpha}^\top\mathbf{d} + \bolds{
\beta}^\top\mathbf{b} \bigr) &=& \exp\Biggl( \sum
_{i,j=1;i\neq j}^n (\alpha_i+
\beta_j) a_{i,j} \Biggr)
\nonumber\\[-8pt]\\[-8pt]\nonumber
&=& \prod
_{i,j=1;i\neq j}^n \exp\bigl((\alpha_i+
\beta_j) a_{i,j} \bigr),
\end{eqnarray}
which implies that the $n(n-1)$ random variables $a_{i,j}$, $i \neq j$
are mutually independent and $Z(\bolds{\alpha}, \bolds{\beta})$ can be
expressed as
%
\begin{equation}
\label{z-alpha-beta-con} \quad Z(\bolds{\alpha}, \bolds{\beta}) = \sum
_{i\neq j} Z_1(\alpha_i+\beta
_j):= \sum_{i\neq j} \log\int
_\Omega\exp\bigl( (\alpha_i+\beta
_j) a_{i,j} \bigr) \nu(da_{i,j}).
\end{equation}
Since an out-edge from vertex $i$ pointing to $j$ is the in-edge of $j$
coming from $i$, it is immediate that
\[
\sum_{i=1}^n d_i = \sum
_{j=1}^n b_j.
\]
Moreover, since the sample is just one realization of the random graph,
the density or probability mass function \eqref{Eq:density:whole} is
also the likelihood function. Note that if one transforms $(\bolds
{\alpha
}, \bolds{\beta})$ to $(\bolds{\alpha}-c, \bolds{\beta}+c)$, the
likelihood does not change. Therefore, for identifiability, constraints
on $\bolds{\alpha}$ or $\bolds{\beta}$ are necessary. In this paper, we
choose to set $\beta_n=0$. Other constraints are also possible, for
example, $\sum_i\alpha_i=0$ or $\sum_j\beta_j=0$.
In total, there are $2n-1$ independent parameters and the natural
parameter space becomes
\[
\Theta=\bigl\{ (\alpha_1, \ldots, \alpha_n,
\beta_1, \ldots, \beta_{n-1})^\top\in
R^{2n-1}: Z(\bolds{\alpha}, \bolds{\beta})<\infty\bigr\}.
\]
Note that model \eqref{Eq:density:whole} can serve as the null model
for hypothesis testing, for example, \cite
{HollandLeinhardt1981,FienbergWasserman1981b}, or be used to
reconstruct directed networks and make statistical inference in a
situation in which only the bi-degree sequence is available due to
privacy consideration \cite{HelleringerKohler2007}. Moreover,
many complex directed network models reply on the bi-degree sequences,
indirectly
or directly. Thus, model \eqref{Eq:density:whole} can be used for
preliminary analysis of network data for
choosing suitable statistics in describing network configurations, for
example, \cite{Robins.et.al.2009}.

It is worth to note that the above discussions only consider
independent edges. Exponential random graph models (ERGMs), sometimes
referred as exponential-family random graph models, for example, \cite
{HunterHandcock2006,Schweinberger2011}, can be more general.
If dependent network configurations such as $k$-stars and triangles are
included as sufficient statistics, then edges are not independent and
such models incur ``near-degeneracy'' in the sense of \cite
{Handcock2003}, in which almost all realized graphs essentially either
contain no edges or are complete \cite
{Handcock2003,Schweinberger2011,ChatterjeeDiaconis2013}.
It has been shown in \cite{ChatterjeeDiaconis2013} that most
realizations from many ERGMs look similar to the results of a simple
Erd\H{o}s--R\'{e}nyi model, which implies that many distinct models
have essentially the same MLE, and it was also proved and characterized
in \cite{ChatterjeeDiaconis2013} the degeneracy observed in the ERGM
with the counts of edges and triangles as the exclusively sufficient statistics.
Further, by assuming a finite dimension of the parameter space, it was
shown in \cite{ShaliziRinaldo2013} that the MLE is not consistent in
the ERGM when the sufficient statistics involve $k$-stars, triangles
and motifs of $k$-nodes ($k\ge2$), while it is consistent when edges
are dyadic independent. In view of the model degeneracy and problematic
properties of estimators in the ERGM for dependent network
configurations, we choose not to consider dependent edges in this paper.


For the remainder of the paper, we proceed as follows.
In Section~\ref{section:main}, we first introduce notation and key
technical propositions that will be used in the proofs. We establish
asymptotic results in the cases of binary weights, continuous weights
and discrete weights in Sections~\ref{subsection:binary},~\ref
{subsection:continuous} and~\ref{subsection:discrete}, respectively.
Simulation studies are presented in Section~\ref{section:simulation}.
We further discuss the results in Section~\ref{section:discussion}.
Since the technical proofs in Sections~\ref{subsection:continuous} and
\ref{subsection:discrete} are similar to those in Section~\ref
{subsection:binary}, we show the proofs for the theorems in
Section~\ref{subsection:binary} in the \hyperref[appendix]{Appendix}, while the proofs for
Sections~\ref{subsection:continuous} and~\ref{subsection:discrete},
as well as those for Proposition~\ref{pro:inverse:appro}, Theorem~\ref
{theorem:Newton:converg} and Lemmas~\ref{Lemma-binary-1} and~\ref
{Lemma-binary-2} in Section~\ref{subsection:binary}
are relegated to the Online Supplementary Material
\cite{YanLengZhu2014}.

\section{Main results}\label{sec2}
\label{section:main}

\subsection{Notation and preparations}
\label{subsection:notations}

Let\vspace*{1pt} $\R_+ = (0, \infty)$, $\R_0=[0, \infty)$, $\N= \{1,\break 2,\ldots\}
$, $\N_0 = \{0,1,2 ,\ldots\}$. For a subset $C\subset\R^n$, let
$C^0$ and $\overline{C}$ denote the interior and closure of $C$,
respectively. For a vector $\mathbf{x}=(x_1, \ldots, x_n)^\top\in
R^n$, denote by
$\llVert \mathbf{x}\rrVert _\infty= \max_{1\le i\le n} \llvert
x_i\rrvert $, the $\ell_\infty
$-norm of $\mathbf{x}$. For an $n\times n$ matrix $J=(J_{i,j})$, let
$\llVert J\rrVert _\infty$ denote the matrix norm induced by the $\ell
_\infty
$-norm on vectors in $\R^n$, that is,
\[
\llVert J\rrVert_\infty= \max_{\mathbf{x}\neq0}
\frac{ \llVert J\mathbf{x}\rrVert
_\infty}{\llVert \mathbf{x}\rrVert _\infty} =\max_{1\le i\le n}\sum_{j=1}^n
\llvert J_{i,j}\rrvert.
\]

In order to characterize the Fisher information matrix, we introduce a
class of matrices.
Given two positive numbers $m$ and $M$ with $M \ge m >0$, we say the
$(2n-1)\times(2n-1)$ matrix $V=(v_{i,j})$ belongs to the class
$\mathcal{L}_{n}(m, M)$ if the following holds:
%
\begin{eqnarray}\label{eq1}
m&\le& v_{i,i}-\sum _{j=n+1}^{2n-1} v_{i,j} \le M,\qquad i=1,\ldots,
n-1; \nonumber
\\
v_{n,n}&=&\sum_{j=n+1}^{2n-1} v_{n,j},\nonumber
\\
v_{i,j}&=&0,\qquad i,j=1,\ldots,n, i\neq j,\nonumber
\\
v_{i,j}&=&0,\qquad i,j=n+1, \ldots, 2n-1, i\neq j,
\\
m&\le& v_{i,j}=v_{j,i} \le M,\qquad i=1,\ldots, n, j=n+1,\ldots,
2n-1, j\neq n+i,\hspace*{-20pt}\nonumber
\\
v_{i,n+i}&=&v_{n+i,i}=0,\qquad i=1,\ldots,n-1,\nonumber
\\
v_{i,i}&=& \sum_{k=1}^n v_{k,i}=\sum_{k=1}^n v_{i,k},\qquad i=n+1, \ldots, 2n-1.\nonumber
\end{eqnarray}
%
Clearly, if $V\in\mathcal{L}_{n}(m, M)$, then $V$ is a $(2n-1)\times
(2n-1)$ diagonally dominant, symmetric nonnegative
matrix and $V$ has the following structure:
\[
V= \pmatrix{ V_{11} & V_{12}
\vspace*{3pt}\cr
V_{12}^\top&
V_{22}},
\]
where $V_{11}$ ($n$ by $n$) and $V_{22}$ ($n-1$ by $n-1$) are diagonal
matrices, $V_{12}$ is a nonnegative matrix whose nondiagonal elements
are positive and diagonal elements equal to zero.

Define $v_{2n,i}=v_{i,2n}:= v_{i,i}-\sum_{j=1;j\neq i}^{2n-1} v_{i,j}$
for $i=1,\ldots, 2n-1$ and $v_{2n,2n}=\sum_{i=1}^{2n-1} v_{2n,i}$.
Then $m \le v_{2n,i} \le M$ for $i=1,\ldots, n-1$, $v_{2n,i}=0$ for
$i=n, n+1,\ldots, 2n-1$ and $v_{2n,2n}=\sum_{i=1}^n v_{i, 2n}=\sum
_{i=1}^n v_{2n, i}$. We propose to approximate the inverse of $V$,
$V^{-1}$, by the matrix $S=(s_{i,j})$, which is defined as
\[
s_{i,j}=\cases{ \displaystyle\frac{\delta_{i,j}}{v_{i,i}} + \frac{1}{v_{2n,2n}},
&\quad$i,j=1,
\ldots,n$,
\cr
\displaystyle -\frac{1}{v_{2n,2n}}, &\quad$i=1,\ldots, n$, $j=n+1,\ldots,2n-1$,
\cr
\displaystyle -\frac{1}{v_{2n,2n}}, &\quad$i=n+1,\ldots,2n-1$, $j=1,\ldots,n$,
\cr
\displaystyle\frac{\delta_{i,j}}{v_{i,i}}+\frac{1}{v_{2n,2n}}, &\quad$i,j=n+1,\ldots
, 2n-1$,}
\]
where $\delta_{i,j}=1$ when $i=j$ and $\delta_{i,j}=0$ when $i\neq
j$. Note that $S$ can be rewritten as
\[
S = \pmatrix{ S_{11} & S_{12}
\vspace*{3pt}\cr
S_{12}^\top&
S_{22}},
\]
where $S_{11} =1/v_{2n, 2n} + \operatorname{diag}(1/v_{1,1}, 1/v_{2,2}, \ldots,
1/v_{n,n})$, $S_{12}$ is an $n\times(n-1)$ matrix whose elements are
all equal to $-1/v_{2n, 2n}$, and $S_{22} = 1/v_{2n, 2n} + \operatorname{diag}(1/v_{n+1, n+1}, 1/v_{n+2, n+2}, \ldots, 1/v_{2n-1, 2n-1})$.

To quantify the accuracy of this approximation, we define another
matrix norm $\llVert \cdot\rrVert $ for a matrix $A=(a_{i,j})$ by
$\llVert A\rrVert:=\max_{i,j} \llvert a_{i,j}\rrvert $. Then we have
the following proposition, whose proof
is given in the Online Supplementary Material \cite{YanLengZhu2014}.

\begin{proposition} \label{pro:inverse:appro}
If $V\in\mathcal{L}_n(m, M)$ with $M/m=o(n)$, then for large enough~$n$,
\[
\bigl\llVert V^{-1}-S \bigr\rrVert\le\frac{c_1M^2}{m^3(n-1)^2},
\]
where $c_1$ is a constant that does not depend on $M$, $m$ and $n$.
\end{proposition}

Note that if $M$ and $m$ are bounded constants, then the upper bound of
the above approximation error is on the order of $n^{-2}$, indicating
that $S$ is a high-accuracy approximation to $V^{-1}$.
Further, based on the above proposition, we immediately have the
following lemma.

\begin{lemma}\label{lemma:inverse:bound}
If $V\in\mathcal{L}_n(m, M)$ with $M/m=o(n)$, then for a vector
$\mathbf{x}\in R^{2n-1}$,
\begin{eqnarray*}
\bigl\llVert V^{-1}\mathbf{x} \bigr\rrVert_\infty& \le&
\bigl\llVert\bigl(V^{-1}-S\bigr)\mathbf{x} \bigr\rrVert
_\infty+ \llVert S\mathbf{x}\rrVert_\infty
\\
& \le& \frac{2c_1(2n-1)M^2\llVert \mathbf{x}\rrVert _\infty
}{m^3(n-1)^2}+ \frac
{ \llvert x_{2n}\rrvert }{v_{2n,2n}} + \max_{i=1,\ldots, 2n-1}
\frac{ \llvert x_i\rrvert }{ v_{i,i} },
\end{eqnarray*}
where $x_{2n}:=\sum_{i=1}^n x_i - \sum_{i=n+1}^{2n-1}x_i$.
\end{lemma}

Let $\bolds{\theta}=(\alpha_1, \ldots, \alpha_n, \beta_1, \ldots,
\beta_{n-1})^\top$ and $\mathbf{g}=(d_1, \ldots, d_n, b_1, \ldots,
b_{n-1})^\top$.\break Henceforth, we will use $V$ to denote the Fisher
information matrix of the parameter vector $\bolds{\theta}$ and show
$V\in\mathcal{L}_n(m, M)$. In the next three subsections, we will
analyze three specific choices of the weight set: $\Omega=\{0,1\}$,
$\Omega=\R_0$, $\Omega=\N_0$, respectively. For each case, we
specify the distribution of the edge weights $a_{i,j}$, the natural
parameter space $\Theta$, the likelihood equations, and prove the
existence, uniqueness, consistency and asymptotic normality of the MLE.
We defer the proofs for the results in Section~\ref{subsection:binary}
to the \hyperref[appendix]{Appendix}
and all other proofs for Sections~\ref{subsection:continuous} and~\ref
{subsection:discrete} to the Online Supplementary Material \cite{YanLengZhu2014}.

\subsection{Binary weights}\label{subsection:binary}

In the case of binary weights, that is, $\Omega=\{0,1\}$, $\nu$ is
the counting measure, and $a_{i,j}$, $1\le i\neq j\le n$ are mutually
independent Bernoulli random variables with
\[
P(a_{i,j}=1) = \frac{ e^{\alpha_i + \beta_j} }{ 1 + e^{\alpha_i +
\beta_j} }.
\]
%
The log-partition function $Z(\bolds{\theta})$ is $\sum_{i\neq j}\log(
1 + e^{\alpha_i+\beta_j} )$ and the likelihood equations are
%
\begin{eqnarray}
\label{eq:likelihood-binary} d_i & = & \sum_{k=1, k\neq i}^n
\frac{ e^{\hat{\alpha}_i +
\hat{\beta}_k } }{
1+ e^{\hat{\alpha}_i+\hat{\beta}_k}},\qquad  i=1,\ldots, n,
\nonumber\\[-8pt]\\[-8pt]\nonumber
b_j & = & \sum_{k=1,k\neq j}^n
\frac{ e^{\hat{\alpha}_k +
\hat{\beta}_j }}{
1 + e^{\hat{\alpha}_k + \hat{\beta}_j } },\qquad j=1,\ldots, n-1,
\end{eqnarray}
where $\hat{\bolds{\theta}}=(\hat{\alpha}_1, \ldots, \hat
{\alpha}_n, \hat{\beta}_1, \ldots, \hat{\beta}_{n-1} )^\top$ is
the MLE of $\bolds{\theta}$
and $\hat{\beta}_n=0$. Note that in this case, the likelihood
equations are identical to the moment equations.

We first establish the existence and consistency of $\hat{\bolds
{\theta}}$ by applying Theorem~\ref{theorem:Newton:converg} in the \hyperref[appendix]{Appendix}.
Define a system of functions:
\begin{eqnarray*}
F_i(\bolds{\theta}) & = & d_i - \sum
_{k=1; k \neq i}^n \frac{e^{\alpha
_i+\beta_k}}{1+e^{\alpha_i+\beta_k} },\qquad  i=1,\ldots, n,
\\
F_{n+j}(\bolds{\theta}) & = & b_j - \sum
_{k=1; k\neq j}^n \frac
{e^{\alpha_k+\beta_j}}{1+e^{\alpha_k+\beta_j}},\qquad j=1,\ldots, n,
\\
F(\bolds{\theta}) & = & \bigl(F_1(\bolds{\theta}), \ldots,
F_{2n-1}(\bolds{\theta})\bigr)^\top.
\end{eqnarray*}
Note the solution to the equation $F(\bolds{\theta})=0$ is precisely the
MLE. Then the Jacobian matrix $F'(\bolds{\theta})$ of $F(\bolds{\theta})$
can be calculated as follows. For $i=1,\ldots, n$,
\begin{eqnarray*}
\frac{\partial F_i }{\partial\alpha_l} &=& 0,\qquad l=1,\ldots, n, l\neq i;\qquad
\frac{\partial F_i}{\partial\alpha_i}= - \sum
_{k=1;k \neq i}^n \frac{e^{\alpha_i+\beta_k}}{(1+e^{\alpha_i+\beta_k})^2},
\\
\frac{\partial F_i}{\partial\beta_j} &=& -\frac{e^{\alpha_i+\beta
_j}}{(1+e^{\alpha_i+\beta_j})^2},\qquad j=1,\ldots, n-1, j\neq i;\qquad
\frac{\partial F_i}{\partial\beta_i}=0
\end{eqnarray*}
and for $j=1,\ldots, n-1$,
\begin{eqnarray*}
\frac{\partial F_{n+j} }{\partial\alpha_l} &=& -\frac{e^{\alpha
_l+\beta_j}}{(1+e^{\alpha_l+\beta_j})^2},\qquad l=1,\ldots, n, l\neq j;\qquad
\frac{\partial F_{n+j}}{\partial\alpha_j} =0,
\\
\frac{\partial F_{n+j}}{\partial\beta_j} &=& - \sum_{k=1;k \neq j}^n
\frac{e^{\alpha_k+\beta_j}}{(1+e^{\alpha_k+\beta_j})^2},\qquad \frac
{\partial F_{n+j} }{\partial\beta_l} = 0,\qquad l=1,\ldots, n-1.
\end{eqnarray*}
First, note that since the Jacobian is diagonally dominant with
nonzero diagonals, it is positive definite, implying that the
likelihood function has a unique optimum.\vadjust{\goodbreak} Second, it is not difficult
to verify that $-F'(\bolds{\theta})\in\mathcal{L}_n(m, M)$, thus
Proposition~\ref{pro:inverse:appro} and Theorem~\ref
{theorem:Newton:converg} can be applied. Let $\bolds{\theta}^*$ denote
the true parameter vector.
The constants $K_1$, $K_2$ and $r$ in the upper bounds of Theorem~\ref
{theorem:Newton:converg} are given in the following lemma,
whose proof is given in the Online Supplementary Material \cite{YanLengZhu2014}.

\begin{lemma}\label{Lemma-binary-1}
Take $D=R^{2n-1}$ and $\bolds{\theta}^{(0)}=\bolds{\theta}^*$ in Theorem
\ref{theorem:Newton:converg}.
Assume
%
\begin{eqnarray}\label{assumption:binary:a}
&& \max\Bigl\{ \max_{i=1,\ldots,n}\bigl\llvert
d_i - \E(d_i)\bigr\rrvert, \max_{j=1,\ldots,
n}
\bigl\llvert b_j - \E(b_j) \bigr\rrvert\Bigr\}
\nonumber\\[-8pt]\\[-8pt]\nonumber
&&\qquad \le \sqrt{(n-1)\log(n-1)}.
\end{eqnarray}
Then we can choose the constants $K_1$, $K_2$ and $r$ in Theorem~\ref
{theorem:Newton:converg} as
\[
K_1=n-1,\qquad
K_2=\frac{n-1}{2},\qquad
r\le \frac{(\log
n)^{1/2}}{n^{1/2}} \bigl(c_{11} e^{6\llVert \bolds{\theta}^*\rrVert
_\infty} + c_{12}e^{2\llVert \bolds{\theta}^*\rrVert _\infty}
\bigr),
\]
where $c_{11}$ and $c_{12}$ are constants.
\end{lemma}

The following lemma assures that condition \eqref{assumption:binary:a}
holds with a large probability,
whose proof is again given in the Online Supplementary Material \cite{YanLengZhu2014}.

\begin{lemma}\label{Lemma-binary-2}
With probability at least $1-4n/(n-1)^2$, we have
\[
\max\Bigl\{ \max_i\bigl\llvert d_i -
\E(d_i) \bigr\rrvert, \max_j \bigl\llvert
b_j - \E( b_j)\bigr\rrvert\Bigr\} \le\sqrt{(n-1)
\log(n-1) }.
\]
\end{lemma}

Combining the above two lemmas, we have the result of consistency.

\begin{theorem}\label{Theorem:binary:con}
Assume that $\bolds{\theta}^*\in\R^{2n-1}$ with $\llVert \bolds{\theta
}^*\rrVert
_\infty\le\tau\log n $, where $0<\tau<1/24$ is a constant,
and that $A \sim\P_{\bolds{\theta}^*}$, where $\P_{\bolds{\theta}^*}$ denotes
the probability distribution \eqref{Eq:density:whole} on $A$ under the
parameter $\bolds{\theta}^*$. Then as $n$ goes to infinity,
with probability approaching one, the MLE $\hat{\bolds{\theta}}$ exists
and satisfies
\[
\bigl\llVert\hat{\bolds{\theta}} - \bolds{\theta}^* \bigr\rrVert
_\infty= O_p \biggl( \frac{ (\log n)^{1/2}e^{8\llVert \bolds{\theta
}^*\rrVert _\infty} }{ n^{1/2} }
\biggr)=o_p(1).
\]
Further, if the MLE exists, it is unique.
\end{theorem}

Next, we establish asymptotic normality of $\hat{\bolds{\theta}}$
and outline the main ideas in the following. Let $\ell(\bolds{\theta};
A)=\sum_{i=1}^n \alpha_i d_i + \sum_{j=1}^{n-1}\beta_j b_j -\sum_{i\neq
j}\log(1+e^{\alpha_i+\beta_j})$ denote the log-likelihood
function of the parameter vector $\bolds{\theta}$ given the sample $A$.
Note that $F'(\bolds{\theta})=\partial^2 \ell/\partial\bolds{\theta
}^2$, and $V=-F'(\bolds{\theta})$ is the Fisher information matrix of
the parameter vector $\bolds{\theta}$. Clearly, $\hat{\bolds{\theta
}}$ does not have an explicit expression according to the system of
likelihood equations \eqref{eq:likelihood-binary}. However, if
$\hat{\bolds{\theta}}$ can be approximately represented as a
function of $\mathbf{g}=(d_1, \ldots, d_n, b_1, \ldots,
b_{n-1})^\top$ with an explicit expression, then the central limit
theorem for $\hat{\bolds{\theta}}$ immediately follows by noting
that under certain\vadjust{\goodbreak} regularity conditions\vspace*{-1pt}
\[
\frac{g_i - \E(g_i) }{v_{i,i}^{1/2}} \to N(0,1),\qquad n\to\infty,
\]
where $g_i$ denotes the $i$th element of $\mathbf{g}$. The identity
between the likelihood equations and the moment equations provides such
a possibility. Specifically, if we apply Taylor's expansion to each
component of $\mathbf{g}-\E(\mathbf{g})$, the second-order term in
the expansion is $V(\hat{\bolds{\theta}} - \bolds{\theta})$, which
implies that obtaining an expression of $\hat{\bolds{\theta}} -
\bolds
{\theta}$ crucially depends on the inverse of $V$. Note that
$V=-F'(\bolds{\theta})\in\mathcal{L}_n(m, M)$ according to the
previous calculation. Although $V^{-1}$ does not have a closed form, we
can use $S$ to approximate it and Proposition~\ref{pro:inverse:appro}
establishes an upper bound on the error of this approximation, which is
on the order of $n^{-2}$ if $M$ and $m$ are bounded constants.

Regarding the asymptotic normality of $g_i - \E(g_i)$, we note that
both $d_i=\sum_{k\neq i}a_{i,k}$ and $b_j=\sum_{k\neq j}a_{k,j}$
are sums of $n-1$ independent Bernoulli random variables. By the central
limit theorem for the bounded case in \cite{Loeve1977}, page~289,
we know that $v_{i,i}^{-1/2} (d_i - \E(d_i) )$ and
$v_{n+j,n+j}^{-1/2}(b_j - \E(b_j))$ are asymptotically standard normal
if $v_{i,i}$ diverges. Since $e^x/(1+e^x)^2$ is an increasing function
on $x$ when $x\ge0$ and
a decreasing function when $x\le0$, we have
\[
\frac{(n-1)e^{2\llVert \bolds{\theta}^*\rrVert _\infty
}}{(1+e^{2\llVert \bolds{\theta}^*\rrVert
_\infty})^2}\le v_{i,i} \le\frac{n-1}{4},\qquad  i=1, \ldots,
2n.
\]
%
In all, we have the following proposition.

\begin{proposition}\label{pro:binary:central}
Assume that $A\sim\P_{\bolds{\theta}^*}$.
If $e^{\llVert \bolds{\theta}^*\rrVert _\infty}=o( n^{1/2} )$, then
for any fixed
$k \ge1$, as
$n\to\infty$, the vector consisting of the first $k$ elements of
$S\{\mathbf{g}-\E(\mathbf{g})\}$ is asymptotically multivariate
normal with mean zero
and covariance matrix given by the upper left $k \times k$ block of $S$.
\end{proposition}

The central limit theorem is stated in the following and proved by
establishing a relationship between $\hat{\bolds{\theta}}-\bolds
{\theta}$ and
$S\{\mathbf{g}-\E(\mathbf{g})\}$ (see details in the \hyperref[appendix]{Appendix} and
the Online Supplementary Material \cite{YanLengZhu2014}).

\begin{theorem}\label{Theorem:binary:central}
Assume that $A\sim\P_{\bolds{\theta}^*}$. If $\llVert \bolds{\theta
}^*\rrVert
_\infty\le\tau\log n$, where $\tau\in(0,\break  1/44)$ is a constant,
then\vspace*{1pt} for any fixed $k\ge1$, as $n \to\infty$, the vector consisting
of the first $k$ elements of $(\hat{\bolds{\theta}}-\bolds{\theta
}^*)$ is asymptotically multivariate normal with mean $\mathbf{0}$ and
covariance matrix given by the upper left $k \times k$ block of $S$.
\end{theorem}

\begin{remark}
By Theorem~\ref{Theorem:binary:central}, for any fixed $i$, as
$n\rightarrow\infty$, the convergence rate of $\hat{\theta}_i$ is
$1/v_{i,i}^{1/2}$. Since $(n-1)e^{-2\llVert \bolds{\theta}^*\rrVert
_\infty}/4\le
v_{i,i}\le(n-1)/4$, the rate of convergence is between
$O(n^{-1/2}e^{\llVert
\bolds{\theta}^*\rrVert _\infty})$ and $O(n^{-1/2})$.
\end{remark}

In this subsection, we have presented the main ideas to prove the
consistency and asymptotic normality of the MLE for the case of binary
weights. In the next\vadjust{\goodbreak} two subsections, we apply similar ideas to the
cases of continuous and discrete weights, respectively.

\subsection{Continuous weights}
\label{subsection:continuous}

Another important case of model \eqref{Eq:density:whole} is when the
weight of the edge is continuous. For example, in communication
networks, if an edge denotes the talking time between two people in a
telephone network, then its weight is continuous. In the case of
continuous weights, that is, $\Omega=[0, \infty)$, $\nu$ is the
Borel measure and $a_{i,j}$, $1\le i\neq j\le n$ are mutually
independent exponential random variables with the density
\[
f_{\bolds{\theta}}(a)= \frac{1}{-(\alpha_i+\beta_j)}e^{(\alpha_i +
\beta_j)a},\qquad
\alpha_i + \beta_j < 0,
\]
and the natural parameter space is
\[
\Theta=\{\bolds{\theta}: \alpha_i + \beta_j <0\}.
\]
To follow the tradition that the rate parameters are positive in
exponential families, we take the transformation $\bar{\bolds{\theta
}}=-\bolds{\theta}$, $\bar{\alpha}_i = - \alpha_i$ and $\bar{\beta
}_j=-\beta_j$. The corresponding natural parameter space then becomes
\[
\overline{\Theta}=\{\bar{\bolds{\theta}}: \bar{\alpha}_i + \bar{
\beta}_j > 0\}.
\]
Here, we denote by $\hat{\bolds{\theta}}$ the MLE of $\bar{\bolds
{\theta}}$. The log-partition $Z(\bar{\bolds{\theta}})$ is $\sum_{i\neq
j}\log( \bar{\alpha}_i + \bar{\beta}_j )$ and the
likelihood equations are
%
\begin{eqnarray}
\label{eq:likelihood:continuous} d_i & = & \sum_{k=1;k\neq i}^n
( \hat{\alpha}_i+\hat{\beta}_k )^{-1},\qquad  i=1,
\ldots, n,
\nonumber\\[-8pt]\\[-8pt]\nonumber
b_j & = & \sum_{k=1;k\neq j}^n (
\hat{\alpha}_k + \hat{\beta}_j )^{-1},\qquad j=1,
\ldots, n.
\end{eqnarray}

Similar to Section~\ref{subsection:binary}, we define a system of functions:
\begin{eqnarray*}
F_i(\bar{\bolds{\theta}}) & = & d_i - \sum
_{k \neq i} ( \bar{\alpha}_i + \bar{
\beta}_k )^{-1},\qquad  i=1,\ldots, n,
\\
F_{n+j}(\bar{\bolds{\theta}} ) & = & b_j - \sum
_{k\neq j} ( \bar{\alpha}_k + \bar{
\beta}_j )^{-1},\qquad j=1,\ldots, n-1,
\\
F(\bar{\bolds{\theta}}) & = & \bigl(F_1( \bar{\bolds{\theta}} ),
\ldots, F_{2n-1}( \bar{\bolds{\theta}} ) \bigr)^\top.
\end{eqnarray*}
The solution to the equation $F(\bar{\bolds{\theta}})=0$ is the MLE,
and the Jacobian matrix $F'(\bar{\bolds{\theta}})$ of $F(\bar{\bolds
{\theta}})$ can be calculated as follows. For $i=1,\ldots, n$,
\begin{eqnarray*}
\frac{\partial F_i }{\partial\bar{\alpha}_l} &=& 0,\qquad l=1,\ldots, n, l\neq
i;\qquad \frac{\partial F_i}{\partial\bar{\alpha}_i} = \sum
_{k \neq i} \frac{1}{(\bar{\alpha}_i + \bar{\beta}_k)^2},
\\
\frac{\partial F_i}{\partial\bar{\beta}_j} &=& \frac{1}{(\bar
{\alpha}_i + \bar{\beta}_j)^2},\qquad j=1,\ldots, n-1, j\neq i;\qquad
\frac{\partial F_i }{\partial\bar{\beta}_i } =0,
\end{eqnarray*}
and for $j=1,\ldots, n-1$,
\begin{eqnarray*}
\frac{\partial F_{n+j}}{\partial\bar{\alpha}_l} &=& \frac{1}{ (\bar
{\alpha}_l + \bar{\beta}_j)^2},\qquad l=1,\ldots, n, l\neq j;\qquad
\frac{\partial F_{n+j} }{\partial\bar{\alpha}_j }=0,
\\
\frac{\partial F_{n+j}}{\partial\bar{\beta}_j} &=& \sum_{k \neq j} \frac
{1}{(\bar{\alpha}_j + \bar{\beta}_j)^2};\qquad
\frac{\partial F_{n+j}}{\partial\bar{\beta}_l} = 0,\qquad l=1,\ldots, n-1,
l\neq j.
\end{eqnarray*}
It is not difficult to see that $F'(\bar{\bolds{\theta}}^*) \in
\mathcal{L}_n(m, M)$ such that Proposition~\ref{pro:inverse:appro}
can be applied,
and the constants in the upper bounds of Theorem~\ref
{theorem:Newton:converg} are given in the following lemma.

\begin{lemma}\label{Lemma-continuous-1}
Assume that $\bar{\bolds{\theta}}^*$ satisfies $q_n\le\bar{\alpha}_i^*
+\bar{\beta}_j^*\le Q_n$ for any $1\le i\neq j\le n$ and
%
\begin{equation}
\label{assumption:continuous:a} \max\Bigl\{ \max_{i=1,\ldots,n}\bigl
\llvert
d_i - \E(d_i)\bigr\rrvert, \max_{j=1,\ldots,
n}
\bigl\llvert b_j - \E(b_j)\bigr\rrvert\Bigr\} \le
\sqrt{ \frac{ 8(n-1)\log n }{ \gamma q_n^2
} },
\end{equation}
where $\gamma$ is an absolute constant.
Then we have
\[
r = \bigl\llVert\bigl[F'\bigl(\bar{\bolds{\theta}}^*\bigr)
\bigr]^{-1}F\bigl(\bar{\bolds{\theta}}^*\bigr)\bigr\rrVert
_\infty\le\biggl( \frac{2c_1Q_n^6}{nq_n^4}+\frac{1}{(n-1)q_n^2} \biggr
) \sqrt
{\frac{8(n-1)\log n}{ \gamma q_n^2 } }.
\]
Further, take $\bar{\bolds{\theta}}^{(0)}=\bar{\bolds{\theta}}^*$
and $D = \Omega(\bar{\bolds{\theta}}^*, 2r)$ in Theorem~\ref
{theorem:Newton:converg}, that is, an open ball $\{\bolds{\theta}:
\llVert
\bolds{\theta} - \bar{\bolds{\theta}}^* \rrVert _{\infty} < 2r \}$.
If $q_n-4r>0$, then we can choose $K_1=2(n-1)/(q_n-4r)^3$ and
$K_2=(n-1)/(q_n-4r)^3$.
\end{lemma}

The following lemma assures condition \eqref{assumption:continuous:a}
holds with a large probability.

\begin{lemma}\label{Lemma-continuous-2}
With probability at least $1-4/n$, we have
\[
\max\Bigl\{ \max_i\bigl\llvert d_i -
\E(d_i)\bigr\rrvert, \max_j \bigl\llvert
b_j - \E(b_j)\bigr\rrvert\Bigr\} \le\sqrt{
\frac{ 8(n-1)\log n }{ \gamma q_n^2 } }.
\]
\end{lemma}

Combining the above two lemmas, we have the result of consistency.

\begin{theorem}\label{Theorem:continuous:con}
Assume that $\bar{\bolds{\theta}}^*$ satisfies $q_n\le\bar{\alpha}_i^*
+\bar{\beta}_j^*\le Q_n$
and $A \sim P_{\bar{\bolds{\theta}}^*}$. If $Q_n/q_n=o\{ (n/\log
n)^{1/18}\}$, then as $n$ goes to infinity,
with probability approaching one, the MLE $\hat{\bolds{\theta}}$ exists
and satisfies
\[
\bigl\llVert\hat{\bolds{\theta}} - \bar{\bolds{\theta}}^* \bigr\rrVert
_\infty= O_p\biggl( \frac{ Q_n^9 (\log n)^{1/2} }{ n^{1/2}q_n^9 }
\biggr)=o_p(1).
\]
Further, if the MLE exists, it is unique.
\end{theorem}

Again, note that both $d_i=\sum_{k\neq i}a_{i,k}$ and $b_j=\sum_{k\neq
j} a_{k,j}$ are sums of $n-1$ independent exponential random
variables, and
$V=F'(\bar{\bolds{\theta}}^*)\in\mathcal{L}_n(m, M)$ is the Fisher
information matrix of $\bar{\bolds{\theta}}$. It is not difficult to
show that the third moment of the exponential random variable with rate
parameter $\lambda$ is $6\lambda^{-3}$. Under the assumption of
$0<q_n\le\bar{\alpha}_i^* + \bar{\beta}_j^* \le Q_n$, we have
\[
\frac{ \sum_{j=1;j\neq i}^n
\E(a_{i,j}^3) }{ v_{i,i}^{3/2} } = \frac{ 6\sum_{j=1;j\neq i}^n
(\bar{\alpha}_i^* + \bar{\beta}_j^*)^{-1} }{ v_{i,i}^{1/2} } \le\frac{
6Q_n/q_n}{(n-1)^{1/2}}\qquad\mbox{for }i=1,\ldots,n
\]
and
\[
\frac{ \sum_{i=1;i\neq j}^n
\E(a_{i,j}^3) }{ v_{n+j,n+j}^{3/2} } = \frac{ 6\sum_{i=1;i\neq j}^n
(\bar{\alpha}_i^* + \bar{\beta}_j^*)^{-1} }{ v_{n+j,n+j}^{1/2} } \le
\frac{ 6Q_n/q_n}{(n-1)^{1/2}}\qquad\mbox{for }j=1, \ldots, n.
\]
Note that if $Q_n/q_n=o(n^{1/2})$, the above expression goes to zero.
This implies that the condition for the Lyapunov's central limit
theorem holds.
Therefore, $v_{i,i}^{-1/2} (d_i - \E(d_i))$ is asymptotically standard
normal if $Q_n/q_n=o(n^{1/2})$.
Similarly, $v_{n+j, n+j}^{-1/2}(b_j - \E(b_j))$ is also asymptotically
standard normal under the same condition. Noting that $[S (\mathbf
{g}-\E(\mathbf{g}) )]_i =v_{i,i}^{-1} (g_i - \E(g_i) ) +
v_{2n,2n}^{-1} (b_n -
\E(b_n) )$, we have the following proposition.

\begin{proposition}\label{pro:continuous:cen}
If $Q_n/q_n=o( n^{1/2} )$, then for any fixed $k\ge1$, as $n\to\infty
$, the vector consisting of the first $k$ elements of $S( \mathbf
{g}-\E(\mathbf{g}) ) $ is asymptotically multivariate normal with
mean zero and covariance matrix given by the upper $k\times k$ block of
the matrix $S$.
\end{proposition}

By establishing a relationship between $\hat{\bolds{\theta}}-\bar{\bolds{\theta}}^*$ and $S\{\mathbf{g}-\E(\mathbf{g}) \}$, we have
the central limit theorem
for the MLE $\hat{\bolds{\theta}}$.

\begin{theorem}\label{Theorem:continuous:central}
If $Q_n/q_n=o( n^{1/50}/(\log n)^{1/25})$, then for any fixed $k\ge1$, as
$n \to\infty$, the vector consisting of the first $k$ elements of
$\hat{\bolds{\theta}} - \bar{\bolds{\theta}}^*$ is asymptotically
multivariate normal with mean zero and covariance matrix given by the
upper $k\times k$ block of the matrix $S$.
\end{theorem}

\begin{remark}
By Theorem~\ref{Theorem:continuous:central}, for any fixed $i$, as $n
\rightarrow\infty$, the convergence rate of $\hat{\theta}_i$ is
$1/v_{i,i}^{1/2}$. Since $(n-1)/Q_n^2\le v_{i,i}\le(n-1)/q_n^2$, the
rate of convergence is between $O(n^{-1/2}Q_n)$ and $O(n^{-1/2}q_n)$.
\end{remark}

\subsection{Discrete weights}
\label{subsection:discrete}

In the case of discrete weights, that is, $\Omega=\N_0$, $\nu$ is
the counting measure and $a_{i,j}$, $1\le i\neq j\le n$ are mutually
independent geometric random variables with the probability mass function
\[
P(a_{i,j}=a) = \bigl(1 - e^{(\alpha_i+\beta_j)}\bigr)e^{(\alpha_i +
\beta
_j)a},\qquad a=0,
1, 2,\ldots,
\]
where $\alpha_i + \beta_j < 0$. The natural parameter space is
$\Theta=\{\bolds{\theta}: \alpha_i + \beta_j <0\}$. Again, we take
the transformation $\bar{\bolds{\theta}}=-\bolds{\theta}$, $\bar{\alpha
}_i = - \alpha_i$ and $\bar{\beta}_j=-\beta_j$, and the
corresponding natural parameter space becomes
\[
\overline{\Theta}=\{\bar{\bolds{\theta}}: \bar{\alpha}_i + \bar{\beta
}_j > 0\}.
\]
The log-partition $Z(\bar{\bolds{\theta}})$ is $\sum_{i\neq j}\log(
1- e^{-(\bar{\alpha}_i + \bar{\beta}_j)})$ and the likelihood
equations are
\begin{large}
%
\begin{eqnarray}
\label{eq-likeli-dis} \quad\qquad d_i & = & \sum_{k\neq i}
\frac{ e^{-(\hat{\alpha}_i+\hat{\beta}_k)} }{1 - e^{-(\hat{\alpha
}_i+\hat{\beta}_k)} } =\sum_{k\neq i} \frac{1}{e^{(\hat{\alpha}_i+\hat
{\beta
}_k)} -1},\qquad  i=1,\ldots, n,
\\
b_j & = & \sum_{k\neq j}
\frac{ e^{-(\hat{\alpha}_k+\hat{\beta}_j)} }{1 - e^{-(\hat{\alpha
}_k+\hat{\beta}_j)} } =\sum_{k\neq j} \frac{1}{e^{(\hat{\alpha}_k+\hat
{\beta
}_j)} -1},\qquad j=1,
\ldots, n-1.
\end{eqnarray}
\end{large}

We first establish the existence and consistency of $\hat{\bolds
{\theta}}$ by applying Theorem~\ref{theorem:Newton:converg}. Define a
system of functions:
\begin{eqnarray*}
F_i(\bar{\bolds{\theta}}) & = & d_i - \sum
_{k \neq i} \frac
{1}{e^{(\bar{\alpha}_i+\bar{\beta}_k)} -1},\qquad  i=1,\ldots, n,
\\
F_{n+j}(\bar{\bolds{\theta}}) & = & b_j - \sum
_{k\neq j} \frac
{1}{e^{(\bar{\alpha}_k+\bar{\beta}_j)} -1},\qquad j=1,\ldots, n,
\\
F(\bar{\bolds{\theta}}) & = & \bigl(F_1( \bar{\bolds{\theta}} ),
\ldots, F_{2n-1}(\bar{\bolds{\theta}})\bigr)^\top.
\end{eqnarray*}
The solution to the equation $F(\bar{\bolds{\theta}})=0$ is the MLE,
and the Jacobian matrix $F'(\bar{\bolds{\theta}})$ of $F(\bar{\bolds
{\theta}})$ can be calculated as follows: for $i=1,\ldots, n$,
\begin{eqnarray*}
\frac{\partial F_i }{\partial\bar{\alpha}_l} &=& 0,\qquad l=1,\ldots, n, l\neq
i;\qquad \frac{\partial F_i}{\partial\bar{\alpha}_i} = \sum
_{k=1; k \neq i}^n \frac{e^{(\bar{\alpha}_i+\bar{\beta}_k)}
-1}{(e^{(\bar{\alpha}_i+\bar{\beta}_k)} -1)^2},
\\
\frac{\partial F_i}{\partial\bar{\beta}_j} &=& \frac{e^{(\bar
{\alpha}_i+\bar{\beta}_j)} -1}{
(e^{(\bar{\alpha}_i+\bar{\beta}_j)} -1)^2},\qquad j=1,\ldots, n-1, j\neq i;\qquad
\frac{\partial F_i }{\partial\bar{\beta}_i } =0,
\end{eqnarray*}
and for $j=1,\ldots, n-1$,
\begin{eqnarray*}
\frac{\partial F_{n+j} }{\partial\bar{\alpha}_l} &=& \frac{e^{(\bar
{\alpha}_l+\bar{\beta}_j)} -1}{[e^{(\bar{\alpha}_l+\bar{\beta
}_j)} -1]^2},\qquad l=1,\ldots, n, l\neq j;\qquad
\frac{\partial F_{n+j}
}{\partial\bar{\alpha}_j } =0,
\\
\frac{\partial F_{n+j}}{\partial\bar{\beta}_j } &=& \sum_{k \neq j} \frac
{e^{(\bar{\alpha}_k+\bar{\beta}_j)} -1}{[e^{(\bar{\alpha
}_k+\bar{\beta}_j)} -1]^2};\qquad
\frac{\partial F_{n+j} }{\partial\bar{\beta}_l} = 0,\qquad l=1,\ldots, n-1,
l\neq j.
\end{eqnarray*}
Let $\bar{\bolds{\theta}}^*$ be the true parameter vector. It is not
difficult to see $F'(\bar{\bolds{\theta}}^*)\in\mathcal{L}_n(m, M)$
so that Proposition~\ref{pro:inverse:appro} can be applied. The
constants in the upper bounds of Theorem~\ref{theorem:Newton:converg}
are given in the following lemma.

\begin{lemma}\label{Lemma-discrete-1}
Assume that $\bar{\bolds{\theta}}^*$ satisfies $q_n\le\bar{\alpha
}_i^* + \bar{\beta}_j^* \le Q_n$ for all $i\neq j$, $A\sim\P_{\bar{\bolds{\theta}}^*}$ and
%
\begin{equation}
\label{assumption:discrete:a} \max\Bigl\{ \max_{i=1,\ldots,n}\bigl
\llvert
d_i - \E(d_i)\bigr\rrvert, \max_{j=1,\ldots,
n}
\bigl\llvert b_j - \E(b_j)\bigr\rrvert\Bigr\} \le
\sqrt{ \frac{ 8(n-1)\log n }{ \gamma q_n^2
} },
\end{equation}
where $\gamma$ is an absolute constant.
Then we have
\[
r = \bigl\llVert\bigl[F'\bigl(\bar{\bolds{\theta}}^*\bigr)
\bigr]^{-1}F\bigl(\bar{\bolds{\theta}}^*\bigr)\bigr\rrVert
_\infty\le O \biggl( q_n^{-1}\bigl(
e^{3Q_n}\bigl(1+q_n^{-4}\bigr)+ e^{Q_n}
\bigr)\sqrt{\frac{\log n}{n}} \biggr).
\]
Further, take $\bar{\bolds{\theta}}^{(0)}=\bar{\bolds{\theta}}^*$ and
$D = \Omega(\bar{\bolds{\theta}}^*, 2r)$ in Theorem~\ref
{theorem:Newton:converg}, that is, an open ball $\{\bolds{\theta}:
\llVert \bolds
{\theta} - \bar{\bolds{\theta}}^* \rrVert _{\infty} < 2r \}$.
If $q_n-4r>0$, then we can choose $K_1 =
2(n-1)e^{q_n-4r}(1+e^{q_n-4r})(e^{q_n-4r}-1)^{-2}$ and $K_2
=(n-1)e^{q_n-4r}(1+e^{q_n-4r})(e^{q_n-4r}-1)^{-2}$.
\end{lemma}

The following lemma assures that the condition in the above lemma holds
with a large probability.

\begin{lemma}\label{Lemma-discrete-2}
With probability at least $1-4n/(n-1)^2$, we have
\[
\max\Bigl\{ \max_i\bigl\llvert d_i -
\E(d_i)\bigr\rrvert, \max_j \bigl\llvert
b_j - \E(b_j)\bigr\rrvert\Bigr\} \le\sqrt{
\frac{ 8(n-1)\log n }{ \gamma q_n^2 } }.
\]
\end{lemma}

Combining the above two lemmas, we have the result of consistency.

\begin{theorem}\label{Theorem:discrete:con}
Assume that $\bar{\bolds{\theta}}^*$ satisfies $q_n\le\bar{\alpha
}_i^* + \bar{\beta}_j^* \le Q_n$ for all $i\neq j$ and $A \sim P_{\bar{\bolds{\theta}}^*}$. If $(1+q_n^{-11})e^{6Q_n} = o(n^{1/2}/(\log
n)^{1/2})$ then as $n$ goes to infinity, with probability approaching
one, the MLE $\hat{\bolds{\theta}}$ exists and satisfies
\[
\bigl\llVert\hat{\bolds{\theta}} - \bar{\bolds{\theta}}^* \bigr\rrVert
_\infty= O_p \biggl( e^{3Q_n}\biggl(1+
\frac{1}{q_n^5}\biggr)\sqrt{\frac{\log n}{n}} \biggr)=o_p(1).
\]
Further, if the MLE exists, it is unique.
\end{theorem}

Note that both $d_i=\sum_{j\neq i}a_{i,j}$ and $b_j=\sum_{i\neq
j}a_{i,j}$ are sums of $n-1$ independent geometric random variables.
Also note that $q_n\le\bar{\alpha}_i^* + \bar{\beta}^*_j \le Q_n$
and $V=F'(\bar{\bolds{\theta}}^*)\in\mathcal{L}_n(m, M)$, thus we have
\begin{eqnarray*}
\frac{e^{Q_n}}{(e^{Q_n} - 1)^2 } &\le& v_{i,j}\le\frac{e^{q_n} }{
(e^{q_n} - 1)^2 },\qquad  i=1,\ldots,n,  j=n+1,\ldots, 2n, j\neq n+i,
\\
\frac{(n-1)e^{Q_n}}{(e^{Q_n} - 1)^2 } &\le& v_{i,i} \le\frac
{(n-1)e^{q_n} }{ (e^{q_n} - 1)^2 },\qquad  i=1, \ldots,
2n.
\end{eqnarray*}
Using the moment-generating function of the geometric distribution, it
is not difficult to verify that
\[
\E\bigl(a_{i,j}^3\bigr) = \frac{1-p_{i,j}}{p_{i,j}}+
\frac
{6(1-p_{i,j})}{p_{i,j}^2}+\frac{6(1-p_{i,j})^2}{p_{i,j}^3},
\]
where $p_{i,j}=1-e^{-(\bar{\alpha}_i^*+\bar{\beta}^*_j)}$. By
simple calculations, we also have
\[
\E\bigl(a_{i,j}^3\bigr) = v_{i,j}\biggl(6 +
\frac{ e^{\bar{\alpha}_i^*+\bar{\beta
}^*_j}-1}{e^{\bar{\alpha}_i^*+\bar{\beta}^*_j}} + \frac{ 6}{
e^{\bar{\alpha}_i^*+\bar{\beta}^*_j} -1 }\biggr).
\]
It then follows
\[
\frac{\sum_{j\neq i} \E(a_{i,j}^3) }{ v_{i,i}^{3/2} } \le\frac{ 7
+ 6(e^{q_n} -1 )^{-1} }{v_{i,i}^{1/2}} \le\frac{ [7 + 6(e^{q_n} -1
)^{-1}](e^{Q_n}-1) }{ n^{1/2} e^{Q_n/2}}.
\]
Note that if $e^{Q_n/2}/q_n = o(n^{1/2})$, the above expression goes to
zero, which implies that the condition for the Lyapunov's central limit
theorem holds. Therefore, for $i=1,\ldots, n$, $v_{i,i}^{-1/2} (d_i
- \E(d_i))$ is asymptotically standard normal if $e^{Q_n/2}/q_n = o(
n^{1/2} )$. Similarly, for $i=1,\ldots, n$, $v_{n+i, n+i}^{-1/2} (b_i
- \E(b_i) )$ is also asymptotically standard normal if $e^{Q_n/2}/q_n
= o( n^{1/2} )$. Therefore, we have the following proposition.

\begin{proposition}\label{pro:discrete:cen}
If $e^{Q_n/2}/q_n=o( n^{1/2} )$, then for any fixed $k\ge1$, as $n\to
\infty$, the vector consisting of the first $k$ elements of
$S \{\mathbf{g} - \E(\mathbf{g} )\}$ is asymptotically multivariate
normal with mean zero and covariance matrix given by the upper $k\times
k$ block of the matrix $S$.
\end{proposition}

The central limit theorem for the MLE $\hat{\bolds{\theta}}$ is
stated as follows.

\begin{theorem}\label{Theorem:discrete:central}
If $e^{9Q_n}(1+q_n^{-15}) =o\{n^{1/2}/\log n \}$, then for any fixed
$k\ge1$, as $n \to\infty$, the vector consisting of the first $k$
elements of $\hat{\bolds{\theta}} - \bar{\bolds{\theta}^*}$ is
asymptotically multivariate normal with mean zero and covariance matrix
given by the upper $k\times k$ block of the matrix $S$.
\end{theorem}

\begin{remark}
By Theorem~\ref{Theorem:discrete:central}, for any fixed $i$, as
$n\rightarrow\infty$, the
convergence rate of $\hat{\theta}_i$ is $1/v_{i,i}^{1/2}$. Since
$(n-1)e^{Q_n}(e^{Q_n} - 1)^{-2} \le v_{i,i}\le(n-1)e^{q_n}(e^{q_n} -
1)^{-2}$, the rate of convergence is
between $O(n^{-1/2}e^{Q_n/2})$ and $O(n^{-1/2}e^{q_n/2})$.
\end{remark}

\emph{Comparison to} \cite{Haberman1977,Haberman1981}.
It is worth to note that \cite{Haberman1977} proved uniform
consistency and asymptotic normality of the MLE in the Rasch model for
item response theory under the assumption that all unknown parameters
are bounded by a constant.
Further, Haberman (\cite{Haberman1981}, page~60) wrote that
``Since Holland and
Leinhardt's $p_1$ model is an example of an exponential response
model\ldots'' and ``The situation in the
Holland--Leinhardt\vadjust{\goodbreak} model is very
similar, for their model under $\rho=0$ is mathematically equivalent
to the incomplete Rasch model with $g=h$ and $X_{ii}$ unobserved.''
Consequently, it was claimed that the method in \cite{Haberman1977}
can be extended to derive the consistency and asymptotic normality of
the MLE of the $p_1$ model without reciprocity, but a formal proof was
not given. However, these conclusions seem premature due to the
following reasons. First, in an item response experiment, a~total of
$g$ people give answers ($0$ or $1$) to a total of $h$ items. The
outcomes of the experiment naturally form a bipartite undirected graph,
for example, \cite{BollaElbanna2014}, while model \eqref
{Eq:density:whole} is directed.
Second, each vertex in the Rasch model is only assigned one parameter
measuring either the out-degree effect for people or
the in-degree effect for items, while there are two parameters in model
\eqref{Eq:density:whole}, one for the in-degree and the other for the
out-degree, for each vertex simultaneously. Therefore, model \eqref
{Eq:density:whole} cannot be simply viewed as an equivalent Rasch model.
We also note that \cite{Fischer1981} pointed out that the Rasch model
can be considered as the Bradley--Terry model \cite{bradleyterry52}
for incomplete paired comparisons, for which 
\cite{simons1999} proved uniform consistency and asymptotic normality
for the MLE with a diverging number of parameters.
Third, in contrast to the proofs in \cite{Haberman1977}, our methods
utilize an approximate inverse of the Fisher information matrix,
requiring no upper bound on the parameters, while the methods in \cite
{Haberman1977} were based on the classical exponential family theory of
\cite{BarndorffNielsen1973,Berk1972}.
Therefore, we conjecture that the methods in \cite{Haberman1977}
cannot be extended to study the model in~\eqref{Eq:density:whole}.

\section{Simulation studies}\label{section:simulation}

In this section, we evaluate the asymptotic results for model \eqref
{Eq:density:whole} through numerical simulations.
The settings of parameter values take a linear form. Specifically, for
the case with binary weights,
we set $\alpha_{i+1}^* = (n-1-i)L/(n-1)$ for $i=0, \ldots, n-1$;
for the case with discrete weights, we set $\bar{\alpha}_{i+1}^*
=0.2+ (n-1-i)L/(n-1)$ for $i=0, \ldots, n-1$.
In both cases, we considered four different values for $L$, $L=0$,
$\log(\log n)$, $(\log n)^{1/2}$ and $\log n$, respectively.
For the case with continuous weights, we set $\bar{\alpha}_{i+1}^*=1
+ (n-1-i) L /(n-1)$ for $i=0,\ldots, n-1$ and also four values of $L$
are considered: $L=0$, $\log(\log(n))$, $\log(n)$ and $n^{1/2}$. For
the parameter values of $\bar{\bolds{\beta}}$,
let $\bar{\beta}_i^*=\bar{\alpha}_i^*$, $i=1, \ldots, n-1$ for
simplicity and $\bar{\beta}_n^*=0$ by default.

\begin{figure}

\includegraphics{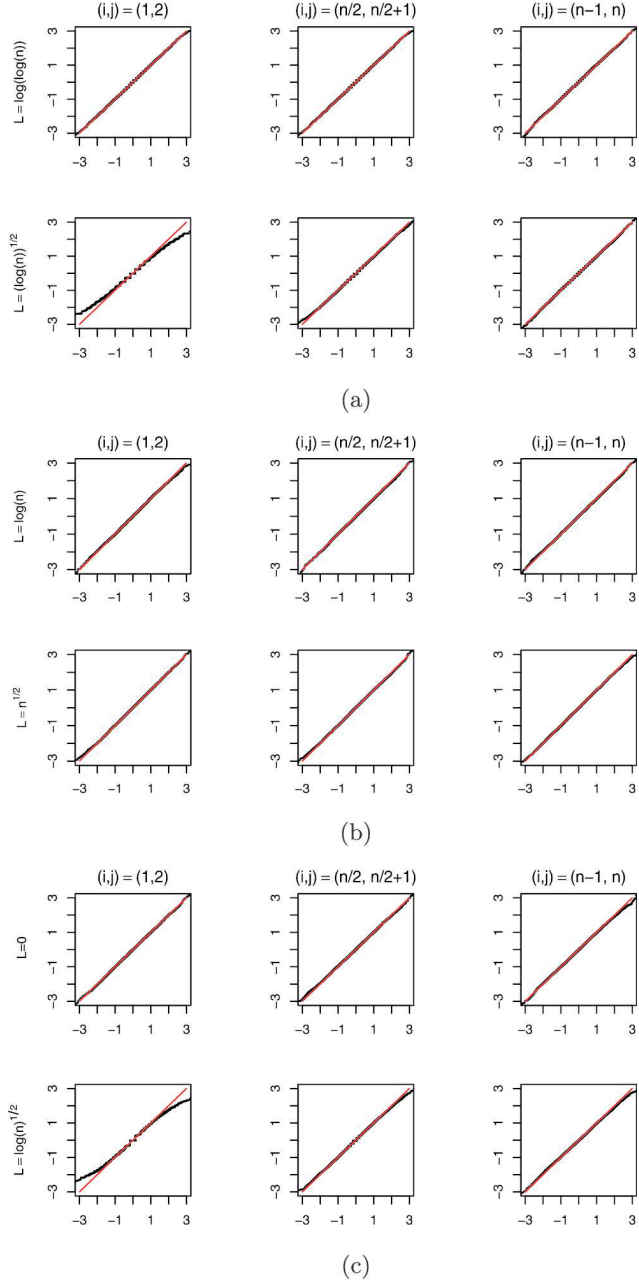}

\caption{The QQ-plots of $\hat{\xi}_{i, j}$ ($n=200$).
\textup{(a)}~Binary weights.
\textup{(b)}~Continuous weights.
\textup{(c)}~Infinite discrete weights.}\label{figure-qq}
\end{figure}

Note that by Theorems~\ref{Theorem:binary:central},~\ref
{Theorem:continuous:central} and~\ref{Theorem:discrete:central}, $\hat
{\xi}_{i,j} = [\hat{\alpha}_i-\hat{\alpha}_j-(\bar{\alpha
}_i^*-\bar{\alpha}_j^*)]/(1/\hat{v}_{i,i}+1/\hat{v}_{j,j})^{1/2}$,
$\hat{\zeta}_{i,j} = (\hat{\alpha}_i+\hat{\beta}_j-\bar{\alpha
}_i^*-\bar{\beta}_j^*)/(1/\hat{v}_{i,i}+1/\hat{v}_{n+j,n+j})^{1/2}$,
and $\hat{\eta}_{i,j} = [\hat{\beta}_i-\hat{\beta}_j-(\bar{\beta
}_i^*-\bar{\beta}_j^*)]/(1/\hat{v}_{n+i,n+i}+1/\hat
{v}_{n+j,n+j})^{1/2}$ are all asymptotically
distributed as standard normal random variables, where $\hat{v}_{i,i}$
is the estimate of $v_{i,i}$ by replacing $\bar{\bolds{\theta}^*}$ with
$\hat{\bolds{\theta}}$. Therefore, we assess the asymptotic
normality of $\hat{\xi}_{i,j}$, $\hat{\zeta}_{i,j}$ and $\hat{\eta
}_{i,j}$ using the quantile--quantile (QQ) plot. Further, we also record
the coverage probability of the 95\% confidence interval, the length of
the confidence interval and the frequency that the MLE does not exist.
The results for $\hat{\xi}_{i,j}$, $\hat{\zeta}_{i,j}$ and $\hat{\eta
}_{i,j}$ are similar, thus only the results of $\hat{\xi
}_{i,j}$ are reported.
Each simulation is repeated 10,000 times.\vadjust{\goodbreak}

We consider two values for $n$, $n=100$ and $200$ and find that the
QQ-plots for them are similar.
Therefore, we only show the QQ-plots when $n=200$ in Figure~\ref
{figure-qq} to save space. In this figure,
the horizontal and vertical axes are the theoretical and empirical
quantiles, respectively,
and the straight lines correspond to the reference line $y=x$.
In Figure~\ref{figure-qq}(b), we can\vspace*{2pt} see that when the weights are
continuous and $L=\log n$ and $n^{1/2}$, the empirical quantiles
coincide with the theoretical ones very well [the QQ-plots when $L=0$
and $\log(\log n)$ are similar to those of $L=\log n$ and not shown].
On the other hand, for binary and discrete weights,
when $L=0$ and $\log(\log n)$, the empirical quantiles agree well with
the theoretical ones
while there are notable deviations when $L=(\log n)^{1/2}$; again, to
save space, the QQ-plots for $L=0$ in the case of binary weights and
for $L=\log(\log n)$ in the case of discrete weights are not shown.
When $L=\log n$, the MLE did not exist in all repetitions (see
Table~\ref{Table:continuous:dis}, thus the corresponding QQ-plot could
not be shown).

\begin{table}
\tabcolsep=0pt
\caption{The reported values are the coverage frequency ($\times$100\%)
for $\alpha_i-\alpha_j$ for a pair $(i,j)$/the length of the
confidence interval/the frequency ($\times$100\%) that the MLE did
not exist}\label{Table:continuous:dis}
\begin{tabular*}{\tablewidth}{@{\extracolsep{\fill}}@{}lccccc@{}}
\hline
$\bolds{n}$ & $\bolds{(i,j)}$ & $\bolds{L=0}$ & $\bolds{L=\log( \log n)}$ & $\bolds{L=(\log(n))^{1/2}}$ & $\bolds{L=\log(n)}$\\
\hline
\multicolumn{6}{@{}c@{}}{\textit{Binary weights}}\\
100 & $(1,2)$ & 94.81/0.57/0 & 95.63/0.10/0.30 & 98.60/1.46/15.86 & NA/NA/100 \\
& $(50,51)$ & 94.78/0.57/0 & 95.18/0.76/0.30 & 95.41/0.93/15.86 & NA/NA/100 \\
& $(99,100)$ & 94.87/0.57/0 & 95.02/0.63/0.30 & 94.97/0.68/15.86 & NA/NA/100
\\[3pt]
200 & $(1,2)$ & 95.35/0.40/0 & 95.50/0.75/0 & 98.13/1.10/1.02 & NA/NA/100 \\
& $(100,101)$ & 95.03/0.40/0 & 95.08/0.55/0 & 95.23/0.68/1.02 & NA/NA/100\\
& $(199,200)$ & 95.28/0.40/0 & 95.32/0.45/0 & 95.26/0.48/1.02 & NA/NA/100
\\[6pt]
\multicolumn{6}{@{}c@{}}{\textit{Continuous weights}}\\
100 & $(1,2)$ & 95.46/1.12/0 & 95.32/2.37/0 & 95.55/4.82/0 & 95.16/9.09/0\\
& $(50,51)$ & 95.28/1.12/0 & 95.44/1.93/0 & 95.71/3.48/0 & 95.51/6.13/0 \\
& $(99,100)$ & 95.38/1.12/0 & 95.63/1.50/0 & 95.81/2.07/0 & 95.72/2.83/0
\\[3pt]
200 & $(1,2)$ & 95.25/0.79/0 & 95.04/1.74/0 & 95.42/3.78/0 & 95.01/8.71/0\\
& $(100,101)$ & 95.10/0.79/0 & 95.21/1.41/0 & 95.31/2.68/0 & 95.39/5.73/0\\
& $(199,200)$ & 95.53/0.79/0 & 95.62/1.07/0 & 95.40/1.52/0 & 95.21/2.28/0
\\[6pt]
\multicolumn{6}{@{}c@{}}{\textit{Discrete weights}} \\
100 & $(1,2)$ & 95.22/0.23/0 & 96.83/1.98/0.54 & 99.72/3.29/56.83 & NA/NA/100 \\
& $(50,51)$ & 95.72/0.23/0 & 95.93/1.15/0.54 & 96.18/1.66/56.83 & NA/NA/100 \\
& $(99,100)$ & 95.49/0.23/0 & 95.73/0.52/0.54 & 95.63/0.61/56.83 & NA/NA/100
\\[3pt]
200 & $(1,2)$ & 95.08/0.16/0 & 96.02/1.51/0 & 98.26/2.56/12.63 & NA/NA/100 \\
& $(100,101)$ & 95.31/0.16/0 & 95.55/0.87/0 & 95.43/1.23/12.63 & NA/NA/100 \\
& $(199,200)$ & 95.28/0.16/0 & 95.54/0.38/0 & 95.31/0.44/12.63 & NA/NA/100 \\
\hline
\end{tabular*}
\end{table}

Table~\ref{Table:continuous:dis} reports the coverage probability of
the 95\% confidence interval for \mbox{$\alpha_i - \alpha_j$}, the length of
the confidence interval, and the frequency that the MLE did not exist.
As we can see, the length of the confidence interval increases as $L$
increases and decreases as $n$ increases, which qualitatively agree
with the theory. In the case of continuous weights, the coverage
frequencies are all close to the nominal level, while in the case of
binary and discrete weights, when $L=(\log n)^{1/2}$ (conditions in
Theorem~\ref{Theorem:discrete:central} no longer hold), the MLE often
does not exist and the coverage frequencies for the $(1, 2)$ pair are
higher than the nominal level; when $L=\log n$, the MLE did not exist
in any of the repetitions.

\section{Summary and discussion}\label{section:discussion}

In this paper, we have derived the uniform consistency and asymptotic
normality of MLEs in the directed ERGM with the bi-degree sequence as
the sufficient statistics; the edge weights are allowed to be binary,
continuous or infinitely discrete and the number of vertices goes to
infinity. In this class of models, a remarkable characterization is
that the Fisher information matrix of the parameter vector is
symmetric, nonnegative and diagonally dominant such that an
approximately explicit expression of the MLE can be obtained.

In the case of discrete weights, only binary and infinitely countable
values have been considered. In the finite discrete case, we may assume
$a_{i,j}$ takes values in the set $\Omega=\{0, 1, \ldots, q-1\}$,
where $q$
is a fixed constant. By \eqref{Eq:density:whole}, it can be shown that
the probability mass function of $a_{i,j}$ is of the form
\[
\P( a_{i,j} = a ) = \frac{ 1 - e^{ - (\alpha_i + \beta_j) }}{
1- e^{-(\alpha_i + \beta_j)q} }\times e^{ -(\alpha_i + \beta_j)a
},\qquad a=0,
\ldots, q-1,
\]
and the likelihood equations become
%
\begin{eqnarray*}
d_i & = & \sum_{j\neq i}
\frac{ 1 - e^{ - (\alpha_i + \beta_j) }}{
1- e^{-(\alpha_i + \beta_j)q} }\sum_{k=0}^{q-1}e^{ -k(\alpha_i +
\beta_j) },
\\
b_j & = & \sum_{i\neq j} \biggl(
\frac{ 1}{ e^{\hat{\alpha}_i +
\hat{\beta}_j } - 1} - \frac{ q}{ e^{(\hat{\alpha}_i +
\hat{\beta}_j)q} - 1}\biggr).
\end{eqnarray*}
%
It can be shown that the Fisher information matrix of $\bolds{\theta}$
is also in the class of matrices $\mathcal{L}_n(m, M)$ under certain
conditions. Therefore, except for some more complex calculations in
contrast with the binary case, there is no extra difficulty to show
that the conditions of Theorem~\ref{Theorem:binary:con} hold, and the
consistency and asymptotic normality of the MLE in the finite discrete
case can also be established.

It is worth noting that the conditions imposed on $q_n$ and $Q_n$ may
not be best possible. In particular, the conditions guaranteeing the
asymptotic normality seem stronger than those guaranteeing the
consistency. For example, in the case of continuous weights, the
consistency requires $Q_n/q_n=(n/\log n)^{1/18}$, while the asymptotic
normality requires $Q_n/q_n=n^{1/50}/\break (\log n)^{1/25}$. Simulation
studies suggest that the conditions on $q_n$ and $Q_n$ might be relaxed.
We will investigate this in future studies and note that the asymptotic
behavior of the MLE depends not only on $q_n$ and $Q_n$, but also on
the configuration of the parameters.

Regarding the $p_1$ model by \cite{HollandLeinhardt1981}, which is
related to model \eqref{Eq:density:whole},
one of the key features of the $p_1$ model is to measure the
dyad-dependent reciprocation by the reciprocity parameter $\rho$. In
the $p_1$ model, there is also another
parameter ($\lambda$) that measures the density of edges, and the
sufficient statistic of the density parameter $\lambda$ is a linear
combination of the in-degrees of vertices and the out-degrees of
vertices. Specifically, the item $\lambda\sum_{i\neq j} a_{i,j} +
\sum_i \alpha_i d_i + \sum_j \beta_j b_j$ in\vspace*{1pt} the $p_1$ model can be
rewritten as $\sum_i (\alpha_i + \lambda+\beta_n)d_i + \sum_j(\beta_j -
\beta_n)b_j$. Therefore, when there is no reciprocity
parameter $\rho$, by taking the transformation of parameters $\tilde
{\alpha}_i = \alpha_i + \lambda+\beta_n$ and $\tilde{\beta}_j =
\beta_j - \beta_n$, we obtain the model \eqref{Eq:density:whole}.
If the reciprocity parameter is incorporated into model \eqref
{Eq:density:whole}, the induced Fisher information matrix is no longer
diagonally dominant and Proposition~\ref{pro:inverse:appro} cannot be
applied. However, simulation results in \cite{YanLeng2013} indicate
that the MLEs still enjoy the properties of uniform consistency and
asymptotic normality, in which the asymptotic variances of the MLEs are
the corresponding diagonal elements of the inverse of the Fisher
information matrix.
In order to extend the current work to study the reciprocity parameter,
a new approximate matrix to the inverse of the Fisher information
matrix is needed.
We plan to investigate this problem in further work.

Finally, we note that the results in this paper can be potentially used
to test the fit of the $p_1$ model. For example, the issue of testing
the fit of the $p_1$ model has been discussed in several previous work,
including \cite
{HollandLeinhardt1981,FienbergWasserman1981b,PetrovicRinaldoFienberg2010,FienbergPetrovicRinaldo2011},
but mostly in heuristic ways.
In view of the result in this paper that the MLE enjoys good asymptotic
properties in model \eqref{Eq:density:whole}, the conjectures in the
above references on the asymptotic distribution of the likelihood ratio
test for testing the fit of $p_1$ model seem reasonable. For example,
to test $H_0: \rho=0$ against $H_1: \rho\neq0$, the likelihood ratio
test proposed by \cite{HollandLeinhardt1981} is likely well
approximated by the chi-square distribution with one degree of freedom.

\begin{appendix}\label{appendix}
\section*{Appendix: Proofs of theorems}

In this section, we give proofs for the theorems presented in Section~\ref{sec2}.

\subsection{Preliminaries}

We first present 
the interior mapping theorem of the mean parameter space, and establish
the geometric rate of convergence for the Newton iterative algorithm to
solve a system of likelihood equations that will be used in this section.

\subsubsection{Uniqueness of the MLE}
Let $\sigma_\Omega$ be a $\sigma$-algebra over the set of weight
values $\Omega$ and $\nu$ be a canonical $\sigma$-finite probability
measure on $(\Omega, \sigma_\Omega)$.
In this paper, $\nu$ is the Borel measure in the case of continuous
weight and the counting measure in the case of discrete weight.
Denote $\nu^{n(n-1)}$ by the product measure on $\Omega^{n(n-1)}$.
Let $\mathfrak{P}$ be all the probability distributions on $\Omega
^{n\choose2}$ that are absolutely continuous with respective to
$\nu^{n\choose2}$.
Define the mean parameter space $\M$ to be the set of expected degree
vectors tied to $\bolds{\theta}$ from all distributions $\P\in
\mathfrak{P}$:
\[
\M= \{ \E_\P\mathbf{g} \dvtx\P\in\mathfrak{P} \}.
\]
Since a convex combination of probability distributions in $\mathfrak
{P}$ is also a probability distribution in $\mathfrak{P}$, the set $\M
$ is necessarily convex.
If there is no linear combination of the sufficient statistics in an
exponential family distribution that is constant, then
the exponential family distribution is minimal.
It is true for the probability distribution \eqref{Eq:density:whole}.
If the natural parameter space $\Theta$ is open, then $\P$ is
regular. By the general theory for a regular and minimal exponential
family distribution (Theorem~3.3 of~\cite{WainwrightJordan2008}), the
gradient of the log-partition function maps the natural parameter space
$\Theta$ to the interior of the mean parameter space $\M$, and this mapping
\[
\nabla Z: \Theta\to\M^\circ
\]
is bijective.
Note that the solution to $\nabla Z(\bolds{\theta}) = \mathbf{g}$ is
precisely the MLE of $\bolds{\theta}$.
Thus, we have established the following.

\begin{proposition}\label{Prop:RegularMinimal}
Assume $\Theta$ is open. Then there exists a solution $\bolds{\theta}
\in\Theta$ to the MLE equation $\nabla Z(\bolds{\theta}) = \mathbf
{g}$ if and only if $\mathbf{g} \in\M^\circ$, and if such a
solution exists, it is also unique.
\end{proposition}

\subsubsection{Newton iterative theorem}
\label{subsubsection:Newton}
Let $D$ be an open convex subset of $\R^{2n-1}$, $\Omega(\mathbf{x},
r)$ denote the open ball $\{\mathbf{y}\in\R^{2n-1}: \llVert \mathbf
{x}-\mathbf{y}\rrVert _\infty< r \}$ and $\overline{\Omega(\mathbf{x},
r)}$ be its closure, where $\mathbf{x}\in\R^{2n-1}$.
We will use Newton's iterative sequence to prove the existence and
consistency of the MLE. Convergence properties of the Newton's
iterative algorithm have been studied by many mathematicians
\cite
{Kantorovich1948,Ortega1968,OrtegaRheimboldt1970,Tapia1971,Polyak2006}.
For the ad-hoc system of likelihood equations considered in this paper,
we establish a fast geometric rate of convergence for the Newton's
iterative algorithm given in the following theorem, whose proof
is given in Online Supplementary Materials~\cite{YanLengZhu2014}.

\begin{theorem}\label{theorem:Newton:converg}
Define a system of equations
\begin{eqnarray*}
F_i(\bolds{\theta}) &=& d_i - \sum
_{k=1, k\neq i}^n f(\alpha_i +
\beta_k),\qquad  i=1, \ldots, n,
\\
F_{n+j}(\bolds{\theta}) &=& b_j - \sum
_{k=1, k\neq j}^n f(\alpha_k +
\beta_j),\qquad j=1, \ldots, n-1,
\\
F(\bolds{\theta}) &=& \bigl( F_1(\bolds{\theta}), \ldots,
F_n(\bolds{\theta}), F_{n+1}(\bolds{\theta}), \ldots,
F_{2n-1}(\bolds{\theta})\bigr)^\top,
\end{eqnarray*}
where $f(\cdot)$ is a continuous function with the third derivative.
Let $D\subset\R^{2n-1}$ be a convex set and assume for any $\mathbf
{x}, \mathbf{y}, \mathbf{v}\in D$, we have
%
\begin{eqnarray}
\label{Newton-condition-a} \bigl\llVert\bigl[F'(\mathbf{x}) -
F'(\mathbf{y})\bigr]\mathbf{v}\bigr\rrVert_\infty&\le&
K_1 \llVert\mathbf{x} - \mathbf{y}\rrVert_\infty\llVert
\mathbf{v}\rrVert_\infty,
\\
\label{Newton-condition-b} \max_{i=1,\ldots,2n-1} \bigl\llVert F_i'(
\mathbf{x}) - F_i'(\mathbf{y})\bigr\rrVert
_\infty&\le& K_2 \llVert\mathbf{x} - \mathbf{y}\rrVert
_\infty,
\end{eqnarray}
where $F'(\bolds{\theta})$ is the Jacobin matrix of $F$ on $\bolds
{\theta
}$ and $F_i'(\bolds{\theta})$ is the gradient function of $F_i$ on
$\bolds
{\theta}$.
Consider $\bolds{\theta}^{(0)}\in D$ with $\Omega(\bolds{\theta}^{(0)},
2r) \subset D$,
where $r=\llVert [F'(\bolds{\theta}^{(0)})]^{-1}F(\bolds{\theta
}^{(0)}) \rrVert
_\infty$.
For any $\bolds{\theta}\in\Omega(\bolds{\theta}^{(0)}, 2r)$, we assume
%
\begin{equation}
\label{Newton-condition-c} F'(\bolds{\theta}) \in\mathcal{ L}_n(m,
M)\quad\mbox{or}\quad {-}F'(\bolds{\theta}) \in\mathcal{
L}_n(m, M).
\end{equation}
For $k=1, 2, \ldots,$ define the Newton iterates $\bolds{\theta
}^{(k+1)} = \bolds{\theta}^{(k)} - [ F'(\bolds{\theta}^{(k)})]^{-1}
F(\bolds
{\theta}^{(k)})$.
Let
%
\begin{equation}
\label{definition-Newton-rho} \rho= \frac{ c_1(2n-1)M^2K_1}{ 2m^3n^2 }
+ \frac{ K_2 }{ (n-1)m}.
\end{equation}
If $\rho r < 1/2$, then $\bolds{\theta}^{(k)}\in\Omega(\bolds{\theta
}^{(0)}, 2r)$, $k=1, 2, \ldots,$ are well defined and satisfy
%
\begin{equation}
\label{Newton-convergence-rate} \bigl\llVert\bolds{\theta}^{(k+1)} -
\bolds{
\theta}^{(0)} \bigr\rrVert_\infty\le r/(1-\rho r).
\end{equation}
Further, $\lim_{k\to\infty} \bolds{\theta}^{(k)}$ exists and the
limiting point is precisely the solution of $F(\bolds{\theta})=0$
in the range of $\bolds{\theta} \in\Omega(\bolds{\theta}^{(0)}, 2r)$.
\end{theorem}

\subsection{Proofs of Theorems \texorpdfstring{\protect\ref{Theorem:binary:con}}{1} and
\texorpdfstring{\protect\ref{Theorem:binary:central}}{2}}

\subsubsection{Proof of Theorem \texorpdfstring{\protect\ref{Theorem:binary:con}}{1}}
Assume that condition \eqref{assumption:binary:a} holds. Recall the
Newton's iterates $\bolds{\theta}^{(k+1)}=\bolds{\theta}^{(k)} -
[F'(\bolds
{\theta}^{(k)})]^{-1}F(\bolds{\theta}^{(k)})$ with $\bolds{\theta
}^{(0)}=\bolds{\theta}^*$. If $\bolds{\theta} \in\Omega(\bolds{\theta
}^*, 2r)$, then $-F'(\bolds{\theta})\in\mathcal{L}_n(m, M)$ with
\[
M=\frac{1}{4},\qquad m=\frac{ e^{2(\llVert \bolds{\theta}^*\rrVert _\infty
+2r) } }{ (
1 + e^{2(\llVert \bolds{\theta}^*\rrVert _\infty+2r) } )^2 }.
\]
If\vspace*{1pt} $\llVert \bolds{\theta}^*\rrVert _\infty\le\tau\log n$ with the
constant $\tau
$ satisfying $0 < \tau< 1/16$, then as $n\to\infty$,
$ n^{-1/2}(\log n)^{1/2}e^{8\llVert \bolds{\theta}^*\rrVert }\le
n^{-1/2+8\tau
}(\log n)^{1/2} \to0$. By Lemma~\ref{Lemma-binary-1} and condition~\eqref{assumption:binary:a}, for sufficiently small $r$,
\begin{eqnarray*}
\rho r & \le& \biggl[\frac{ c_1(2n-1)M^2(n-1)}{2m^3n^2} + \frac{
(n-1) }{ 2m(n-1) } \biggr]
\\
&&{}\times\frac{(\log n)^{1/2}}{n^{1/2}} \bigl( c_{11} e^{6\llVert
\bolds{\theta}^*\rrVert
_\infty} +
c_{12}e^{2\llVert \bolds{\theta}^*\rrVert _\infty} \bigr)
\\
& \le& O \biggl( \frac{(\log n)^{1/2}e^{12\llVert \bolds{\theta
}^*\rrVert _\infty
}}{ n^{1/2} } \biggr) + O \biggl(\frac{(\log n)^{1/2}e^{8\llVert \bolds
{\theta}^*\rrVert _\infty}}{
n^{1/2} }
\biggr).
\end{eqnarray*}
Therefore, if $\llVert \bolds{\theta}^*\rrVert _\infty\le\tau\log n$,
then $\rho
r\to0$ as $n\to\infty$.
Consequently, by Theorem~\ref{theorem:Newton:converg}, $\lim_{n\to
\infty}\hat{\bolds{\theta}}^{(n)}$ exists. Denote the limit as
$\hat{\bolds{\theta}}$, then it satisfies
\begin{eqnarray*}
&&\bigl\llVert\hat{\bolds{\theta}} - \bolds{\theta}^*\bigr\rrVert
_\infty\le2r = O \biggl(\frac{ (\log n)^{1/2}e^{8\llVert \bolds{\theta
}^*\rrVert _\infty} }{ n^{1/2} } \biggr) = o(1).
\end{eqnarray*}
By Lemma~\ref{Lemma-binary-2}, condition (\ref{assumption:binary:a})
holds with probability approaching one, thus the above inequality also
holds with probability approaching one. The uniqueness of the MLE comes
from Proposition~\ref{Prop:RegularMinimal}.

\subsubsection{Proof of Theorem \texorpdfstring{\protect\ref{Theorem:binary:central}}{2}}
Before proving Theorem~\ref{Theorem:binary:central}, we first establish
two lemmas.

\begin{lemma}\label{central:binary:lemma1}
Let $R = V^{-1}-S$ and $U=\operatorname{Cov}[R \{\mathbf{g} - \E\mathbf{g}\}
]$. Then
%
\begin{equation}
\llVert U\rrVert\le\bigl\llVert V^{-1}-S\bigr\rrVert+
\frac{(1+e^{2\llVert \bolds{\theta}^*\rrVert _\infty
})^4}{4e^{4\llVert \bolds{\theta}^*\rrVert _\infty}(n-1)^2}.
\end{equation}
\end{lemma}
\begin{pf}
Note that
\[
U = WVW^\top=\bigl(V^{-1}-S\bigr) - S (I - V S),
\]
where $I$ is a $(2n-1)\times(2n-1)$ diagonal matrix, and by
inequality (C3) in \cite{YanLengZhu2014}, we have
\[
\bigl\llvert\bigl\{ S(I - VS) \bigr\}_{i,j}\bigr\rrvert= \llvert
w_{i,j}\rrvert\le\frac{3(1+e^{2\llVert \bolds{\theta}^*\rrVert
_\infty})^4}{4e^{4\llVert \bolds{\theta
}^*\rrVert _\infty}(n-1)^2}.
\]
Thus,
\begin{eqnarray*}
\llVert U \rrVert&\le&\bigl\llVert V^{-1} - S \bigr\rrVert+ \bigl
\llVert S ( I_{2n-1}-V S )\bigr\rrVert
\\
&\le&\bigl\llVert V^{-1}-S
\bigr\rrVert+\frac{3(1+e^{2\llVert \bolds{\theta}^*\rrVert _\infty
})^4}{4e^{4\llVert
\bolds{\theta}^*\rrVert _\infty}(n-1)^2}.
\end{eqnarray*}\upqed
\end{pf}

\begin{lemma}\label{central:binary:lemma2}
Assume that the conditions in Theorem~\ref{Theorem:binary:con} hold.
If $\llVert \bolds{\theta}^*\rrVert _\infty\le\tau\log n$ and $\tau<
1/40$, then
for any $i$,
%
\begin{equation}
\hat{\theta}_i- \theta_i^* =
\bigl[V^{-1}\bigl\{\mathbf{g} - \E(\mathbf{g}) \bigr\}
\bigr]_i + o_p\bigl( n^{-1/2} \bigr).
\end{equation}
\end{lemma}

\begin{pf}
By Theorem~\ref{Theorem:binary:con}, we have
\[
\hat{\rho}_n:=\max_{1\le i\le2n-1} \bigl\llvert
\hat{\theta}_i-\theta_i^*\bigr\rrvert
=O_p\large\biggl(\frac{(\log n)^{1/2}e^{8\llVert \bolds{\theta
}\rrVert _\infty} }{
n^{1/2}} \large\biggr).
\]
Let $\hat{\gamma}_{i,j}=\hat{\alpha}_i + \hat{\beta
}_j -\alpha_i-\beta_j$.
By Taylor's expansion, for any $1\le i\neq j\le n$,
\[
\frac{ e^{\hat{\alpha}_i +\hat{\beta}_j } }{ 1+
e^{\hat{\alpha}_i+\hat{\beta}_j} } -\frac{e^{\alpha_i^*+\beta
_j^*}}{1+e^{\alpha_i^*+\beta_j^*}} = \frac{e^{\alpha_i^*+\beta
_j^*}}{(1+e^{\alpha_i^*+\beta
_j^*})^2}\hat{
\gamma}_{i,j} +h_{i,j},
\]
where
\[
h_{i,j}=\frac{e^{\alpha_i^*+\beta_j^*+\phi_{i,j}\hat{\gamma
}_{i,j}}(1-e^{\alpha_i^*+\beta_j^*+\phi_{i,j}\hat{\gamma
}_{i,j}}) }{2(1+e^{\alpha_i^*+\beta_j^*+\phi_{i,j}\hat{\gamma
}_{i,j}})^3}\hat{\gamma}_{i,j}^2,
\]
and $0\le\phi_{i,j}\le1$. By the likelihood equations \eqref
{eq:likelihood-binary}, we have
\[
\mathbf{g} - \E( \mathbf{g} ) = V\bigl(\hat{\bolds{\theta}}-\bolds{
\theta}^*\bigr) + \mathbf{h},
\]
where $\mathbf{h}=(h_1, \ldots, h_{2n-1})^\top$ and,
\begin{eqnarray*}
h_i &=& \sum_{k=1,k\neq i}^n
h_{i,k},\qquad  i=1, \ldots, n,
\\
h_{n+i} &=& \sum_{k=1, k\neq i}^n
h_{k,i},\qquad  i=1, \ldots, n-1.
\end{eqnarray*}
Equivalently,
%
\begin{equation}
\label{Eq-binary-central-theta} \hat{\bolds{\theta}}-\bolds{\theta}^* =
V^{-1}\bigl(
\mathbf{g} - \E( \mathbf{g}) \bigr) + V^{-1}\mathbf{h}.
\end{equation}
Since $\llvert e^{x}(1-e^{x})/(1+e^x)^3\rrvert \le1$, we have
\[
\llvert h_{i,j}\rrvert\le\bigl\llvert\hat{\gamma}_{i,j}^2
\bigr\rrvert/2\le2\hat{\rho}_n^2,\qquad \llvert
h_i\rrvert\le\sum_{j\neq i}\llvert
h_{i,j}\rrvert\le2(n-1)\hat{\rho}_n^2.
\]
Note that
$(S \mathbf{h})_i = h_i/v_{i,i} + (-1)^{1_{\{i>n\}}} h_{2n}/v_{2n,2n}$,
and $(V^{-1} \mathbf{h})_i=(S \mathbf{h})_i+(R \mathbf{h})_i$.
By direct calculations, we have
\[
\bigl\llvert(S \mathbf{h})_i\bigr\rrvert\le\frac{\llvert h_i\rrvert
}{v_{i,i}}+
\frac{ \llvert h_{2n}\rrvert
}{v_{2n,2n}} \le\frac{16\hat{\rho}_n^2(1+e^{2\llVert \bolds{\theta
}^*\rrVert
_\infty})^2}{e^{2\llVert \bolds{\theta}^*\rrVert _\infty}} \le O\large
\biggl( \frac{ e^{20\llVert \bolds{\theta}^*\rrVert _\infty} \log n
}{ n }
\large\biggr),
\]
and by Proposition~\ref{pro:inverse:appro}, we have
\[
\bigl\llvert(R \mathbf{h})_i\bigr\rrvert\le\llVert R \rrVert
_\infty\times\Bigl[(2n-1)\max_i\llvert
h_i\rrvert\Bigr] \le O\biggl( \frac{ e^{22\llVert \bolds{\theta
}^*\rrVert _\infty} \log n }{ n } \biggr).
\]
If $\llVert \bolds{\theta}^*\rrVert _\infty\le\tau\log n$ and $\tau<
1/44$, then
\[
\bigl\llvert\bigl(V^{-1}h\bigr)_i\bigr\rrvert\le\bigl
\llvert(Sh)_i\bigr\rrvert+\bigl\llvert(Rh)_i\bigr
\rrvert=o\bigl(n^{-1/2}\bigr).
\]
This completes the proof.
\end{pf}

\begin{pf*}{Proof of Theorem~\ref{Theorem:binary:central}}
By \eqref{Eq-binary-central-theta}, we have
\[
(\hat{\bolds{\theta}}-\bolds{\theta})_i= \bigl[S\bigl\{
\mathbf{g} - \E( \mathbf{g}) \bigr\}\bigr]_i+ \bigl[R \bigl\{
\mathbf{g} - \E(\mathbf{g}) \bigr\}\bigr]_i + \bigl(V^{-1}
\mathbf{h}\bigr)_i.
\]
By Lemmas~\ref{central:binary:lemma1} and~\ref
{central:binary:lemma2}, if $\llVert \bolds{\theta}^*\rrVert _\infty
\le\tau\log
n$ and $\tau< 1/44$,
then
\[
(\hat{\bolds{\theta}}-\bolds{\theta})_i=\bigl[S \bigl\{
\mathbf{g} - \E( \mathbf{g}) \bigr\}\bigr]_i+o_p
\bigl(n^{-1/2}\bigr).
\]
Therefore, Theorem~\ref{Theorem:binary:central} follows
directly from Proposition~\ref{pro:binary:central}.
\end{pf*}
\end{appendix}

\section*{Acknowledgments}
We thank Runze Li for the role he played as Editor, an Associate Editor
and two referees for their valuable comments and suggestions that have
led to significant improvement of the manuscript.

\begin{supplement}
\stitle{Supplement to ``Asymptotics in directed exponential random graph models with an
increasing bi-degree sequence.''}
\slink[doi]{10.1214/15-AOS1343SUPP} 
\sdatatype{.pdf}
\sfilename{aos1343\_supp.pdf}
\sdescription{The supplemental\vadjust{\goodbreak} material contains proofs for
the lemmas in Section~\ref{subsection:binary},
the theorems and lemmas
in Sections~\ref{subsection:continuous} and~\ref
{subsection:discrete}, Proposition~\ref{pro:inverse:appro} and Theorem
\ref{theorem:Newton:converg}.}
\end{supplement}



%

\printaddresses
\end{document}